\documentclass[12pt]{article}

\newenvironment{proof}{\noindent
\textit{Proof.}}{\hfill{{$\square$}}\\}

\usepackage{latexsym}
\usepackage{amsmath}
\usepackage{amssymb}
\usepackage{amscd}
\usepackage{multicol}

\newtheorem{theorem}{\bf Theorem}[section]
\newtheorem{corollary}[theorem]{\bf Corollary}
\newtheorem{proposition}[theorem]{\bf Proposition}
\newtheorem{lemma}[theorem]{\bf Lemma}
\newtheorem{definition}[theorem]{\bf Definition}

\def\cR{\mathcal{R}}

\def\t{\mathfrak{t}}

\def\t{\frak{t}}

\def\b{\frak{b}}
\def\u{\frak{u}}

\def\bZ{\mathord{\Bbb Z}}
\def\bR{\mathord{\Bbb R}}
\def\bN{\mathord{\Bbb N}}
\def\bC{\mathord{\Bbb C}}
\def\r1{\sqrt{-1}}

\def\cR{{\mathcal R}}
\def\ed{d_1=m_1d_{l+1},\ldots ,d_l=m_ld_{l+1}}
\def\er{r_1=m_1r_{l+1},\ldots ,r_{l-1}=m_{l-1}r_{l+1}}
\def\cS{{\mathcal S}}
\def\cH{{\mathcal H}}
\def\cD{{\mathcal D}}
\def\Wa{W_{\rm aff}}
\title{Quantum cohomology of the infinite dimensional generalized
flag manifolds}

\author{ Augustin-Liviu  Mare
}
\date{}

\begin{document}

\maketitle
%{\bf Running head:} Infinite dimensional flag manifolds

%{\bf Proofs should be sent to:} Augustin-Liviu  Mare, Department
%of Mathematics,

%University of Toronto,  Toronto, Ontario M5S 3G3,
 %Canada

%e-mail: {\tt amare@math.toronto.edu}

%\newpage

\begin{abstract}
Consider the infinite dimensional flag manifold $LK/T$
corresponding to the simple Lie group $K$ of rank $l$ and with
maximal torus $T$. We show that, for $K$ of type $A$, $B$ or $C$,
if we endow the space $H^*(LK/T)\otimes \bR[q_1,\ldots,q_{l+1}]$
(where $q_1,\ldots,q_{l+1}$ are multiplicative variables) with an
$\bR[\{q_j\}]$-bilinear product  satisfying some simple properties
analogous to the quantum product on $QH^*(K/T)$, then the
isomorphism type of the resulting ring is determined by the
integrals of motion of a certain {\it periodic} Toda lattice
system, in exactly the same way as the  isomorphism type of
$QH^*(K/T)$ is determined by the integrals of motion of the {\it
non-periodic} Toda lattice (see Kim \cite{kim}).  This is an
infinite dimensional extension of the main result of \cite{l2} and
at the same time a generalization of \cite{mt}.

\end{abstract}

\section{Introduction}
Let $K$ be a compact, connected, simply connected Lie group,
 $T\subset K$ a maximal torus and $LK$   the
  group of loops in $K$ (i.e. smooth maps from the circle $S^1$ to
$K$). The evaluation of a loop at the point $1$ defines a
topologically trivial bundle $p:LK/T \rightarrow K/T$ of fiber
$LK/K=\Omega(K)$, where $\Omega(K)$ is the space of based loops in
$K$. Consequently $LK/T$ is homeomorphic to $K/T\times \Omega(K)$.
Hence a simple description in terms of generators and relations of
the cohomology ring with real coefficients
  of $LK/T$  can be deduced from the fact --- which is
a direct consequence of a classical result of Serre ---  that the
loop space $\Omega(K)$ has the rational homotopy type of a direct
product of spaces of the type $\Omega(S^m)$, for certain $m\geq 2$
(recall that $H^*(\Omega(S^m))=\bR[x]$, where $\deg x = m-1$) and
from
 Borel's presentation of $H^*(K/T)$.
More precisely,  if $\t$ denotes the Lie algebra of $T$,
we can identify
$$H_2(K/T,\bZ)=\pi_2(K/T)=\pi_1(T)=\ker (\exp :\t\to T)$$
the latter being just the coroot lattice in $\t$. This leads to
the identification of the basis of $H^2(K/T)$ consisting of the
degree 2 Schubert classes with the fundamental
 weights $\lambda_1, \ldots ,\lambda_l \in \t^* $
associated to a simple root system. We have
$$H^*(K/T)=\bR[\lambda_1,\ldots,\lambda_l]/\langle
u_k(\lambda_1,\ldots,\lambda_l)~|~ 1\le k \le l\rangle$$ where
$u_k(\lambda_1,\ldots,\lambda_l)$, $1\leq k \leq l$ are the
fundamental
 homogeneous generators of  the subring of all $W$-invariant polynomials in
$S(\t^*)=\bR[\lambda_1,\ldots,\lambda_l]$. We consider
$p^*(\lambda_1), \ldots, p^*(\lambda_l)\in H^2(LK/T)$, where
$p^*:H^*(K/T) \to H^*(LK/T)$ is the (injective) map induced by $p$
at the level of the cohomology rings. If  $\cH$ denotes
$H^*(LK/T)$,  we deduce that
\begin{equation}\cH= H^*(LK/T)=(\bR[p^*(\lambda_1), \ldots, p^*(\lambda_l)]\otimes
H^*(\Omega (K))))/\langle u_k(p^*(\lambda_1), \ldots,
p^*(\lambda_l))~|~ 1\le k\le l\rangle
\end{equation}
where $H^*(\Omega(K))$ is just a polynomial ring (see above) whose generators
are not involved in the ideal of relations.

Fix a system of simple roots $\alpha_1,\ldots,\alpha_l \in \t^*$
corresponding to the root system of $K$ and let $-\alpha_{l+1}$ be
the highest root (the convention concerning the minus sign comes
from Goodman and Wallach's definition \cite{gw1} of the periodic
Toda lattice; see also section 2 of our paper).
 Then  $\alpha_1^{\vee},\ldots,\alpha_l^{\vee} \in \t$
are a system of simple roots of the coroot system. The strictly
positive integers $m_1,\ldots,m_l$ determined by
$$-\alpha_{l+1}^{\vee} = m_1\alpha_1^{\vee} + \ldots + m_l\alpha_l^{\vee}$$
will be important objects in our paper. Consider the formal
multiplicative variables $q_1,\ldots,q_{l+1}$ and the space
$\cH\otimes \bR[q_1,\ldots ,q_{l+1}]$ which consists of all
expressions of the type
$$a=\sum_{d=(d_1,\ldots,d_{l+1})\geq 0}a_dq^{d},$$
where $q^d$ denotes $q_1^{d_1}\ldots q_{l+1}^{d_{l+1}}$. We would
like to stress that for any $a\in \cH\otimes \bR[q_1,\ldots
,q_{l+1}]$, we denote by $a_d\in \cH$ the coefficient of $q^d$.

The main result of our paper is:
\begin{theorem}\label{main}  Suppose that the Lie group $K$ is simple of type
$A$, $B$ or $C$. Let $\bullet$ be an $\bR[\{q_i\}]$-linear product
on $\cH\otimes \bR[\{q_i\}]$ with the following properties:
\begin{itemize}
\item[(i)] $\bullet$ preserves the grading induced by the
usual grading of $H^*(LK/T)$ combined with $\deg q_j=4$, $1\leq j
\leq l+1$
\item[(ii)] $\bullet$ is a deformation of the usual product, in the sense
that if we formally replace all $q_j$ by $0$, we obtain the usual
product on $\cH$
\item[(iii)] $\bullet$ is commutative
\item[(iv)] $\bullet$ is associative
\item[(v)] $p^*(\lambda_i)\bullet p^*(\lambda_j) =p^*(\lambda_i)
p^*(\lambda_j)+\delta_{ij}q_j+m_im_jq_{l+1}$, $1\le i,j\le l$
\item[(vi)] $(d_i-m_id_{l+1})(p^*(\lambda_j)
\bullet a)_d=(d_j-m_jd_{l+1})(p^*(\lambda_i) \bullet a)_d$, for
 any $d=(d_1,\ldots,d_{l+1})\geq 0$, $a\in \cH$, $1\leq i,j\leq l$
\item[(vii)] For $i\in\{1,\ldots, l\}$, $d=(d_1,\ldots, d_{l+1})\neq 0$,
if the coefficient $(p^*(\lambda_i)\bullet a)_d$ of $q^d$ in
$p^*(\lambda_i)\bullet a$ is nonzero, then $d_i-m_id_{l+1}\ne 0$.
\end{itemize}
Then the ring $\cH\otimes \bR[\{q_j\}]$ is generated by
$p^*(\lambda_1),\ldots ,p^*(\lambda_l)$ and $H^*(\Omega(K)),$
subject to the relations
\begin{equation}\label{relation}F_k(-\langle \alpha_1^{\vee},\alpha_1^{\vee}\rangle q_1,\ldots,
-\langle \alpha_{l+1}^{\vee},\alpha_{l+1}^{\vee}\rangle q_{l+1},
p^*(\lambda_1)\bullet,\ldots,p^*(\lambda_l)\bullet)= 0
\end{equation} for $1\leq k\leq l,$ where $F_k$ are the integrals of
motion of the periodic Toda lattice (see Theorem 1.2 below)
corresponding to the coroot system  of $K$.
\end{theorem}

 The Toda lattice
we are referring to in our theorem is the Hamiltonian system which
consists of the standard symplectic manifold $(\bR^{2l},
\sum_{i=1}^l dr_i\wedge ds_i)$ with the Hamiltonian function
\begin{equation}\label{ham}E=\sum_{i,j=1}^l \langle
\alpha_i^{\vee}, \alpha_j^{\vee}\rangle r_i r_j +\sum_{j=1}^{l+1}
e^{2s_j}, \end{equation} where by definition
\begin{equation}\label{s}s_{l+1}=-m_1s_1 -\ldots -m_ls_l.\end{equation} The following
result, proved by Goodman and Wallach in \cite{gw1} and
\cite{gw3}, gives details concerning the complete integrability of
this system (for more details, see section 2):

\begin{theorem} {\rm (see \cite{gw1}, \cite{gw3})} There exist $l$
functionally independent functions $$E=F_1, F_2, \ldots, F_l
:\bR^{2l}\to \bR,$$ such that for every $k\in\{1,\ldots,l\}$ we have:
\begin{itemize}
\item[(i)] $F_k$ is a homogeneous polynomial in variables $e^{2s_1},\ldots , e^{2s_l},
e^{2s_{l+1}}, r_1,\ldots, r_l,$
\item[(ii)] $\{F_k,E\}=0$,
\item[(iii)] $F_k( 0, \ldots, 0, \lambda_1,\ldots, \lambda_l)=
u_k( \lambda_1,\ldots, \lambda_l)$, where $u_k$ is the $k$-th
fundamental $W$-invariant polynomial (see above).
\end{itemize}
\end{theorem}

{\bf Remarks}

1. So far there exists  no rigorous definition of a product
 $\bullet$ on $\cH\otimes \bR[\{q_j\}]$ with the properties
(i)-(vii) stated in Theorem 1.1. We expect that Gromov-Witten
invariants and a quantum product on $\cH\otimes \bR[\{q_i\}]$
which satisfies those properties can be defined by using the same
arguments as in the case of the finite dimensional flag manifold
$K/T$ (see for instance \cite{fp}), and the following two facts:
\begin{itemize}
\item[a)] the natural complex structure of $LK/T$
(see section 4)
\item[b)] the fact that --- as in the finite dimensional situation
--- the ``first Chern class" of the tangent bundle of $LK/T$ (in the
sense defined by Freed \cite{fr}), equals
$$2\sum_{j=1}^{l+1}\tilde{\lambda}_j,$$
where $\tilde{\lambda}_1,\ldots,\tilde{\lambda}_{l+1}$ denote the
degree 2 Schubert classes  in $LK/T$ (as in the finite dimensional
situation, we identify the corresponding Schubert varieties with
the fundamental weights of the corresponding affine Lie algebra);
this seems to indicate that the
 moduli space of pointed holomorphic maps of given multidegree $d$
 from $\bC P^1$ to $LK/T$ is finite dimensional, of  dimension
 $$4\sum_{j=1}^{l+1}d_j.$$
\end{itemize}

2. Readers who are familiar with the quantum cohomology of $K/T$
will certainly find the requirements (i)-(iv) from Theorem 1.1
quite natural. A few explanations concerning (v), (vi) and (vii)
are needed, though. A crucial relation  is (see Proposition
\ref{pstarl}):
$$p^*(\lambda_i)=\tilde{\lambda_i}-m_i\tilde{\lambda}_{l+1},$$
$1\leq i \leq l$. Hence (v) would follow from
$$\tilde{\lambda}_i\bullet \tilde{\lambda}_j=
 \tilde{\lambda}_{i}\tilde{\lambda}_{j}+
\delta_{ij}q_{j},\quad 1\leq i,j\le l+1$$ and one may expect that
a quantum product on $H^*(LK/T)$ satisfies this equation (for a
simple proof of the same property in the finite dimensional case,
see \cite{kim} or \cite{l1}). In order to explain why should (vi)
be satisfied, too, let us consider the coefficient
$(p^*(\lambda_i)\bullet a)_d\in H^*(LK/T)$ of $q^d$ in
$p^*(\lambda_i)\bullet a$ and $\beta$  an arbitrary element in
$H_*(LK/T)$,  and look at
\begin{equation}\label{a}\langle (p^*(\lambda_i)\bullet a)_d,\beta\rangle
= \langle \tilde{\lambda}_i\bullet a,\beta\rangle
-m_i \langle \tilde{\lambda}_{l+1}\bullet a,\beta\rangle.
\end{equation}
One would expect  the ``Gromov-Witten invariants"
$$\langle \tilde{\lambda}_k\bullet a,\beta\rangle
=\langle \tilde{\lambda}_k| a|\beta\rangle_d$$
to satisfy  the ``divisor property"
$$\langle \tilde{\lambda}_k| a|\beta \rangle_d=
d_k\langle a,\beta\rangle_d ,\quad 1\leq k \leq l+1$$ so that the
right hand side of (\ref{a}) is
\begin{equation}\label{dsubi}(d_i-m_id_{l+1})\langle a,\beta\rangle_d \end{equation} and
property (vi) becomes clear. Property (vii) follows immediately
from (\ref{a}) and (\ref{dsubi}) by taking $\beta\in H_*(LK/T)$
with the property that the left hand side of (\ref{a}) is
non-zero.

3. Even though   the functionally independent integrals of motion
$F_1,\ldots,F_l$ exist {\it for any simple Lie group $K$}, (note
that Theorem 1.2 contains no restriction on $K$), the methods of
our paper enable us to prove the main result only for $K$ of type
$A$, $B$ or $C$ (see Theorem 1.1). More precisely, our proof
relies on the relationship between the polynomials $F_k$ and the
commutator of a certain Laplacian in the universal enveloping
algebra of the corresponding $ax+b$-algebra: it is only for $K$ of
one of the types mentioned before that the  structure of this
commutator has been determined  by Goodman and Wallach in
\cite{gw1} and \cite{gw3} (for more details, see section 2). In
fact, this is the only piece of information which is lacking:
whenever we have elements of the universal enveloping algebra
which commute with the Laplacian,
 Theorem \ref{polynomial} produces relations from  them.

4. Theorem 1.1 is an extension of the theorem of Kim \cite{kim} to
the affine case. The point of view we are adopting here (i.e. an
axiomatic approach to the isomorphism type of the quantum
cohomology ring) is the same as in \cite{l2}. At the same time,
Theorem 1.1  is a generalization of the main result of Guest and
Otofuji \cite{mt}, where the situation of the  periodic flag
manifold $Fl^{(n)}=L(SU(n))/T$ (i.e. when $K$ is of type $A$) has
been considered.

{\bf Acknowledgments:} I would like to thank Martin Guest and
Takashi Otofuji for several discussions on the topics contained in
this paper. I am grateful to Lisa Jeffrey for a careful reading of
the manuscript and suggesting several improvements. I would also
like to thank the referee for several helpful suggestions.

\section{Periodic Toda lattices according to Goodman and Wallach}

The goal of this section is to present some details about the
mechanical system of periodic Toda lattice type, as defined by
Goodman and Wallach \cite{gw1}.

\begin{definition} {\rm (see \cite{gw1})}
 The $(ax+b)$-algebra
corresponding to the extension by $\alpha_{l+1}^{\vee}$
of the Dynkin diagram corresponding to the
coroot system of $K$ is
the Lie algebra  $$(\b=\t^* \oplus \u, [ \ ,\ ]),$$ where the Lie
bracket $[\ , \ ]$ is defined
by:
\begin{itemize}
\item $\t^*$ and $\u$ are abelian
\item $\u$ has a basis $X_1,\ldots,X_{l+1}$ such that
\begin{equation}\label{comm}[\lambda,X_j]=\lambda(\alpha_j^{\vee})X_j,
\quad \lambda\in \t^*,
 1\leq j\leq l+1.\end{equation}
\end{itemize}
\end{definition}
The set $S(\b)$ of polynomial functions on $\b^*$ becomes a
Poisson algebra and by (\ref{comm}) we have
$$\{\lambda_{i_1},\lambda_{i_2}\}=0,
\quad \{\lambda_i,X_j\}=\lambda_i(\alpha_j^{\vee})X_j, \quad
\{X_{j_1},X_{j_2}\}=0,$$ for any $1\leq i_1,i_2,i\leq l$, $1\leq
j,j_1,j_2 \leq l+1$.

On the other hand, one can easily see that the Poisson bracket of
functions on the standard symplectic manifold $(\bR^{2l},
\sum_{i=1}^l dr_i\wedge ds_i)$ satisfies
$$\{r_{i_1},r_{i_2}\}=0, \quad \{r_i, e^{s_j}\}=
\lambda_i(\alpha_j^{\vee})e^{s_j}, \quad \{e^{s_{j_1}},
e^{s_{j_2}}\}=0,$$ for $1\leq i_1,i_2,i\leq l$, $1\leq j,j_1,j_2
\leq l+1$ (for $j=l+1$, we have used  the assumption that
 $s_{l+1}=-\sum_{i=1}^lm_is_i$ (see (\ref{s})) and the fact that
$\lambda_i(\alpha_{l+1}^{\vee}) =-m_i$). Consequently the map
$$X_j \mapsto e^{s_j}, \quad \lambda_i \mapsto r_i, $$
$1\leq i \leq l$ is a  homomorphism of Poisson algebras from
$S(\b)$ to the Poisson subalgebra $\bR[e^{s_1},\ldots ,
e^{s_{l+1}}, r_1, \ldots, r_l]$ of $C^{\infty}(\bR^{2l})$. In this
way, integrals of motion of the Hamiltonian system determined by
(\ref{ham}) can be obtained from elements of the space $S(\b)^{\{,
\}}$ consisting of all elements of $S(\b)$ which $\{, \}$-commute
with the polynomial function on $\b^*$ given by
\begin{equation}\label{function} \sum_{i,j=1}^l \langle \alpha_i^{\vee},
\alpha_j^{\vee}\rangle \lambda_i \lambda_j +\sum_{j=1}^{l+1} X_j^2.
\end{equation}

Now let us consider the universal enveloping algebra $$U(\b)=T
(\b)/\langle x\otimes y-y\otimes x -[x,y], x,y \in \b\rangle$$ and
the  isomorphism
$$J:S(\b) \to U(\b)$$
induced by the symmetrization map followed by the canonical
projection (see [7, Corollary E, \S 17.3]). Since $\t^*$ and
$\u$ are abelian, the element of $S(\b)$ described by
(\ref{function}) is mapped by $J$ to $$\Omega:= \sum_{i,j=1}^l
\langle \alpha_i^{\vee}, \alpha_j^{\vee}\rangle \lambda_i
\lambda_j +\sum_{j=1}^{l+1} X_j^2,$$ the right hand side being
regarded this time as an element of $U(\b)$.

 The complete integrability of the Toda lattice follows from the
 following two theorems of Goodman and Wallach:

\begin{theorem}\label{thmgw1} {\rm (see [4, Theorem 6.4])} If $K$ is
of classical type, then the Poisson bracket
commutator $S(\b)^{\{,\}}$ is mapped by $J$ isomorphically onto
the space $U(\b)^{[,]}$ of all $f\in U(\b)$ with the property that
$[f,\Omega]=0$.
\end{theorem}

\begin{theorem}\label{thmgw2} {\rm (see [5, Theorem 3.3])}
Suppose that $K$ is of type $A$, $B$ or $C$. Let
 $\mu: U(\b)\to U(\t^*)=S(\t^*)$ be the map
induced by the natural Lie algebra homomorphism $\b\to \t^*$. Then
$$\mu(U(\b)^{[,]})\subset S(\t^*)^W.$$ Moreover, there exist
$\Omega=\Omega_1, \ldots, \Omega_l \in U(\b)$, such that for each
$k\in\{1,\ldots ,l\}$ we have:
\begin{itemize}
\item[(i)] $[\Omega_k, \Omega]=0,$
\item[(ii)] $\mu(\Omega_k)=u_k$ and $\deg \Omega_k = \deg u_k$, where $u_k$ is the $k$-th
fundamental generator of $S(\t^*)^W$.
\end{itemize}
Each $\Omega_k$ is contained in the subring of
$U(\b)$ which is spanned
by elements of the form $X^{2I}\lambda^J$.
\end{theorem}

Note  that $U(\b)^{[,]}$ is generated as an algebra by the
elements $\Omega_1,\ldots,\Omega_l$ described in the previous
theorem, plus $X_1^{m_1}\ldots X_l^{m_l}X_{{l+1}}$.

\noindent {\bf Remark.} The integrals of motion of the Toda
lattice  mentioned in
 Theorem 1.2
 are the polynomials $F_k$
 obtained from $J^{-1}(\Omega_k)$ by the transformations
(\ref{comm}).

\section{The degree 2  cohomology modules of $K/T$ and  $LK/T$}

As in the introduction, $K$ is a compact, connected,
simply connected and simple Lie
group and $T\subset K$ a maximal torus.
 The goal of this section is to point out  a relation (see Proposition \ref{pstarl} below)
 between  the
cohomology of the complex flag manifold $K/T$ and the cohomology of
 the  infinite
dimensional version    of it,  $LK/T$.

Consider again a simple root system $\{\alpha_1,\ldots, \alpha_l\}\subset\t^*$
and
 $\{\alpha_1^{\vee}, \ldots,
 \alpha_l^{\vee}\}\subset \t$ the corresponding simple coroot system.
As pointed out in the introduction, there is a natural isomorphism
$\phi$  between $H_2(K/T,\bZ)$
 and the coroot lattice in
$\t$,  the latter  being the same as the integral lattice. More
precisely, $\phi$ is the composition of the Hurewicz isomorphism
$H_2(K/T, \bZ)
 \simeq \pi_2 (K/T)$
with the boundary map $\pi _2(K/T) \rightarrow \pi _1 (T)$ from
the long exact sequence of the bundle $T\rightarrow K \rightarrow
K/T$. Consequently, $H^2(K/T, \bZ)$ can be identified with the
weight lattice, i.e.
 the lattice generated by the
elements $\lambda_i$ of $\t^*$, uniquely determined by
$\lambda_i(\alpha_j^{\vee})=\delta_{ij}$, $1\leq i,j\leq l$.

Similar considerations can be made for $LK/T$, but only after
obtaining a special
presentation of it, as follows:
Take first a central extension
$$S^1 \longrightarrow  \tilde{L}K\stackrel{\pi}{\longrightarrow} LK,$$
in the sense of [12, Chapter 4].
If $T\subset LK$ is the set of constant loops in $T$,
then we get the central
extension
$$S^1 \longrightarrow  \pi^{-1}(T) \stackrel{\pi}{\longrightarrow} T.$$
But the only central extension of the torus $T$ is the $S^1\times T$, hence
 $\pi ^{-1}(T)=S^1 \times T.$
The space $LK/T$ can be identified with $\tilde{L}(K)/\pi ^{-1}(T)$ via
 the map  $\tilde{l}\pi ^{-1}(T) \mapsto \pi (\tilde{l}) T$.

And now, exactly as before for $K/T$, we take the composition of
 the Hurewicz
isomorphism $H_2(\tilde{L}K/\pi^{-1}(T), \bZ) \simeq \pi_2
(\tilde{L}K/\pi^{-1}(T))$ with the boundary map $\pi
_2(\tilde{L}K/\pi^{-1}(T)) \rightarrow \pi _1 (\pi^{-1}(T))$ from
the long exact sequence of the bundle $\pi^{-1}(T)\rightarrow
\tilde{L}K \rightarrow \tilde{L}K/\pi^{-1}(T)$, and denote the
resulting isomorphism by
\begin{equation}\label{tilde}\tilde{\phi} :H_2(\tilde{L}K/\pi^{-1}(T), \bZ) \rightarrow
 \pi _1(\pi^{-1}(T)).\end{equation}

Again the integral lattice $\pi _1(\pi^{-1}(T))$ is the same as
the coroot lattice. The latter has a basis consisting of the
simple coroots
 $$\tilde{\alpha}_1^{\vee}=({\alpha}_1^{\vee},0), \ldots ,
 \tilde{\alpha}_l^{\vee}=({\alpha}_l^{\vee},0),
\tilde{\alpha}_{l+1}^{\vee}=({\alpha}_{l+1}^{\vee},-
\frac{1}{2}\langle \alpha_{l+1}^{\vee},  \alpha_{l+1}^{\vee}\rangle),$$
inside the Lie algebra $\t+\bR$ of $\pi^{-1}(T)$,
where $\alpha_{l+1}$ is the highest root (for more details, see [12, Chapter 4]).
So  the isomorphism
$\tilde{\phi} :H_2(\tilde{L}K/\pi^{-1}(T), \bZ) \rightarrow
  \pi _1(\pi^{-1}(T))$
assigns to any coroot $\alpha^{\vee}$ the loop
 ${\rm exp}(t\alpha^{\vee})$, $t\in [0,1]$.
Consequently, $H^2(LK/T)$ can be identified with  the weight lattice, i.e.
 the  lattice generated by the
elements $\tilde{\lambda}_i$ of $(\t+\bR )^*$, uniquely determined by
 $\tilde{\lambda}_i(\tilde{\alpha}_j^{\vee})=\delta_{ij}$,
$1\leq i,j \leq l+1$.

We are interested in the relationship between
$$H_2(K/T, \bR)=\t, \qquad H^2(K/T, \bR)=\t^*$$
on the one hand and
$$H_2(LK/T, \bR)=\t +\bR , \qquad  \ H^2(LK/T)=(\t +\bR )^*$$
on the other. Consider first the inclusion map $I:K/T\rightarrow
LK/T$, which induces  the maps $I_*$ and $I^*$ by functoriality.
One can easily
 see that:
\begin{itemize}
\item[(i)] $I_*$ is the inclusion map of $\t$ into $\t+\bR$;
\item[(ii)]
$\tilde{\lambda}_i|_{\t}=\lambda_i$, $1\leq i \leq l$ and
 $\tilde{\lambda}_{l+1}|_{\t}=0$,
hence $I^*:H^2(LK/T)\rightarrow H^2(K/T)$ maps $\tilde{\lambda}_i$
 to ${\lambda}_i$,
$1\leq i \leq l$ and $\tilde{\lambda}_{l+1}$ to zero.
\end{itemize}
The main result of the section concerns the imbedding of $H^2(K/T)$
into $H^2(LK/T)$
induced by the map
$p:LK/T \rightarrow K/T,$ $$p(\gamma T):=\gamma(1)T.$$
\begin{proposition}\label{pstarl}
We have that
$$p^*(\lambda_i)=\tilde{\lambda}_i-m_i\tilde{\lambda}_{l+1},\quad 1\leq i \leq l, $$
 where $m_i$ are given by $-\alpha_{l+1}^{\vee}=\sum_{i=1}^l m_i
\alpha_{i}^{\vee}.$
\end{proposition}
\begin{proof} From $p\circ I ={\rm id}$ follows that $I^*\circ p^*={\rm id}$, hence
 $p^*(\lambda_i)$ must be of the form
$\tilde{\lambda}_i+k \tilde{\lambda}_{l+1}$,
$k\in \bZ$. In turn, $k$ can be obtained as follows:
$$\begin{array}{lll}k&=p^*(\lambda _i)( \tilde{\alpha}^{\vee}_{l+1})=
\lambda_i(p_*(({\alpha}_{l+1}^{\vee},-
\frac{1}{2}\langle \alpha_{l+1}^{\vee},  \alpha_{l+1}^{\vee}\rangle)))\\{}&=
\lambda_i(-{\alpha}^{\vee}_{l+1})=-m_i,\end{array}$$
 where we have used the following  property of $p$:
\begin{equation}\label{pstar}
 p_*((0,\frac{1}{2}\langle \alpha_{l+1}^{\vee},  \alpha_{l+1}^{\vee}\rangle))
=0. \end{equation}

In order to prove (\ref{pstar}), we consider the map
$P:\tilde{L}K\rightarrow K$, $$P(\tilde{l}):=\pi(\tilde{l})(1),$$
where $\pi$ is the central extension (see above).
 The pair $(P,p)$ is a morphism between the bundles
$(\tilde{L}K, \tilde{L}K/\pi^{-1}T)$ and $(K, K/T)$. By the functoriality of
 the maps $\phi$ and
$\tilde{\phi}$, we have the following commutative diagram:
$$\pi_1(\pi^{-1}(T)) \stackrel{P_{\#}}{\longrightarrow}   \pi_1(T)$$
$$\tilde{\phi}\uparrow \  \  \  \  \  \  \  \  \ \ \ \ \ \  \uparrow \phi$$
$$H_2(LK/T)  \stackrel{p_{*}}{\longrightarrow}   H_2(K/T)$$
The loop  $$\exp(t (0,\frac{1}{2}\langle \alpha_{l+1}^{\vee},
\alpha_{l+1}^{\vee}\rangle)),
 \quad t\in [0,1]$$ in $\pi^{-1}(T)$ corresponds to the coroot
 $(0,\frac{1}{2}\langle \alpha_{l+1}^{\vee},  \alpha_{l+1}^{\vee}\rangle)$
 in $\t+\bR$.
This loop is obviously mapped to the constant loop $e$ by $\pi$,
 hence by $P$. So
(\ref{pstar}) is true.
\end{proof}

\section{The Schubert basis of $H^*(LK/T)$}

The goal of this section is to describe  the Schubert varieties in
$LK/T$, by following the construction outlined in [12, section
8.7].

Let $G$ be the complexification of $K$ and $B_0\subset G$ a Borel
subgroup. Consider the following subgroups of $LG$:
\begin{itemize}
\item[a)] $L^+G$, the group of loops in $G$ which extend holomorphically
 to the
 unit disc;
\item[b)] $B^+$, the subgroup of $L^+G$ consisting of loops $\gamma$ with
 $\gamma(0)\in B_0$.
\end{itemize}
Then we have the diffeomorphism
$$LK/T\simeq LG/B^+,$$
which induces an action of $B^+$ on $LK/T$.

The {\it affine Weyl group} $\Wa$ is defined as  the semidirect
product of $$W=N_K(T)/T$$ with the lattice in $\t$ generated by
the coroots. One can easily see that $\Wa$ is generated by the
reflections $s_1,\ldots,s_l,s_{l+1}$ in the walls
$\ker\alpha_1,\ldots,\ker\alpha_l$, and $\ker\alpha_{l+1}$ of the
fundamental chamber in $\t$. In particular, to any $w\in \Wa$
corresponds a length $l(w)$, which is the minimal number of
factors in a decomposition of $w$ as a product of $s_j$'s. Since
$K$ is simply connected, the coroot lattice coincides with the
integral lattice $\check{T}={\rm Hom}(S^1, T)$. Consequently we
have that
$$\Wa=(N_K(T)\cdot \check{T})/T\subset LK/T.$$
As in the finite dimensional situation, for any $w\in \Wa$, the $B^+$ orbit
 $$C_w=B^+.w$$ is a complex cell of dimension $l(w)$, called the {\it Bruhat cell}.
 The set
$$\{[\overline{C}_w] ~|~ w\in \Wa\}$$ consisting of the fundamental cycles of the closures
of Bruhat cells is a basis of $H_*(LK/T)$. We denote by
$\sigma_w\in H^*(LK/T)$ the dual of $[\overline{C}_w]$ with
respect to the {\it evaluation pairing} $\langle \ , \ \rangle$.
More precisely, $\sigma_w$ is defined by
$$\langle\sigma_w, [\overline{C}_v]\rangle = \delta_{vw},\quad v,w\in \Wa$$
where $\delta_{vw}$ denotes the Kronecker delta.
The basis $\{\sigma_w ~|~ w\in \Wa\}$ is the {\it Schubert basis} of $H^*(LK/T)$.

The set $\{[\overline{C}_{s_j}] ~|~ 1\leq j \leq l+1\}$ is a basis
of $H_2(LK/T,\bZ)$. The isomorphism $\tilde{\phi}$ defined in
section 3 (see equation (\ref{tilde})) maps any
$[\overline{C}_{s_j}]$ to the coroot $\tilde{\alpha}_j^{\vee}$,
$1\leq j \leq l+1$. This can be proved by using the $SU(2)$
embedding in $\tilde{L}(K)$ associated to the root $\alpha_j$ (see
[12, section 5.2]) and the naturality of $\tilde{\phi}$.
Consequently the  transpose of $\tilde{\phi}$ maps $\sigma_{s_j}$
to the fundamental weight $\tilde{\lambda}_j$, $1\leq j\leq l+1$,
and we  often use this isomorphism as an identification.

\section{Relations in the ring $(\cH\otimes \bR[\{q_j\}], \bullet)$:
proof of the main result}

Let us consider the space $$\bR[\{t_i\}][[\{e^{t_j/2}\}]]$$ of all
formal series of the type
$$g=\sum_{d=(d_1\ldots, d_{l+1})\ge 0}g_de^{td}$$
where $g_d\in \bR[\{t_i\}]$ is a polynomial and $e^{td}$ means
$e^{t_1d_1}\ldots e^{t_{l+1}d_{l+1}}$. This space has an obvious
structure of a commutative and associative algebra with unit. For
any $1\leq i \leq l$, we consider  the derivative-type operator
$\partial_i$ on
 $\bR[\{t_i\}][[\{e^{t_j/2}\}]]$,
which is $\bR$-linear and satisfies
 $$\partial_i(ge^{td/2}):=\frac{\partial}{\partial t_i}(g)e^{td/2}+
\frac{1}{2}(d_i-m_id_{l+1})ge^{td/2},$$ for any
 $g\in \bR[t_1, \ldots , t_l]$ and any $d=(d_1,\ldots, d_{l+1})\ge 0$.
The facts that $\partial_i$ is $\bR$-linear and satisfies  the
Leibniz rule can be easily verified.

Now let us consider the following assignment $\rho$:
$$\rho (\lambda_i):=2\partial_i,  \qquad \rho (X_j)=\frac{2
\sqrt{-1}}{h}
 \sqrt{\langle \alpha_j^{\vee} , \alpha _j^{\vee}\rangle}
e^{\frac{t_j}{2}} ,$$ $ 1\leq i \leq l$, $ 1\leq j \leq l+1$,
where $h$ is a nonzero real parameter. One can easily see that
$\rho$ is actually a representation of the Lie algebra $\b$ on the
space $\bR[\{t_i\}][[\{e^{t_j/2}\}]]$, i.e. it satisfies
$$\rho[\lambda_i,X_j]=[\rho(\lambda_i), \rho(X_j)],\quad 1\le i \le l, 1\le j \le l+1.$$

 Consider $\Omega_k$ defined in Theorem \ref{thmgw2}
and  set
$$D_k=h^{{\rm deg}\Omega_k}\rho(\Omega_k), 1\le k \le l.$$
By the last property mentioned in Theorem \ref{thmgw2},  $D_k$
leaves the subspace $\bR[\{t_i\}][[\{e^{t_j}\}]]$ of
$\bR[\{t_i\}][[\{e^{t_j/2}\}]]$ invariant. Only the action on
 $\bR[\{t_i\}][[\{e^{t_j}\}]]$ will be used later, and on this space,
 $\partial_i$ is defined by
\begin{equation}\label{partial}\partial_i(ge^{td}):=\frac{\partial}{\partial t_i}(g)e^{td}+
(d_i-m_id_{l+1})ge^{td},\qquad 1\leq i\leq l.\end{equation}

Since $F_k$ is homogeneous in variables $e^{s_j}, r_i$, $1\leq
j\le l+1$, $1\leq i \le l$,
 it follows that $\Omega_k$ ---
being  essentially the same as $J(F_k)$ --- has a presentation as
a homogeneous, symmetric polynomial in the variables $X_j,
\lambda_i$. We use the commutation relations (\ref{comm}) in order
to express $\Omega_k$ as a linear combination of  elements of the
form $ X^{2I}\lambda^J$ (see  Theorem \ref{thmgw2}).
 The polynomial expression we
obtain in this way appears as
$$ \Omega_k=F_k(X_i^2, \lambda_i) +f_k(X_i^2, \lambda_i)$$
where $$\deg f_k < \deg F_k.$$
Consequently  $D_k$ appears as a
polynomial expression  $D_k(e^{t_1},\ldots
, e^{t_{l+1}},h\frac{\partial}{\partial t_1}, \ldots ,
h\frac{\partial}{\partial t_l}, h)$, the last ``variable", $h$,
being due to the possible occurrence  of $f_k$.

Let us consider the variables $Q_1,\ldots,
Q_{l+1},\Lambda_1,\ldots, \Lambda_l$. To each polynomial $D\in\bR
[Q_j, \Lambda _i, h]$ we assign the differential operator
$D(e^{t_j}, h\partial_i,h)$ obtained from $D$ as follows: first
write $D$ as a sum of monomials of the type $Q^I\Lambda^Jh^m$, and
then replace $Q_j\mapsto e^{t_j}, \Lambda_i \mapsto h\partial_i$,
$1\le j\le l+1$, $1\le i\le l$. For instance, the differential
operator
$$\frac{1}{4}D_1=\frac{1}{4}h^2\rho(\Omega) = \sum_{i,j=1}^l\langle \alpha_i^{\vee},
\alpha_j^{\vee}\rangle h^2\partial_i\partial_j -
\sum_{j=1}^{l+1}\langle \alpha_j^{\vee}, \alpha_j^{\vee}\rangle
e^{t_j}$$ arises as $H(e^{t_j}, h\partial_i,h)$, where
$$H=\sum_{i,j=1}^l\langle \alpha_i^{\vee}, \alpha_j^{\vee}\rangle
\Lambda_i\Lambda_j - \sum_{j=1}^{l+1}\langle \alpha_j^{\vee},
\alpha_j^{\vee}\rangle Q_j.$$ The main result of this section is:

\begin{theorem}\label{polynomial}
Assume  that the polynomial
 $D\in\bR [Q_1, \ldots ,Q_{l+1}, \Lambda _1, \ldots , \Lambda_{l}, h]$
 satisfies
\begin{equation}\label{deg}{\deg}_{\{Q_j^{1/2},\Lambda_i\}}(D)\leq 2\sum_{i=1}^lm_i
.\end{equation}
If
\begin{itemize}
\item[(a)] $[D(e^{t_j}, h\partial_i,h), H(e^{t_j}, h\partial_i,h)]=0$
for any $h\neq 0$,
\item[(b)] the polynomial $D(0, \ldots ,
 0,\Lambda _1, \ldots , \Lambda_l,h)$
does not depend on $h$,
\item[(c)]
$D(0, \ldots , 0,\lambda_1, \ldots , \lambda_l,0)\in S(\t^*)^W$,
\end{itemize}
 then the relation $$D(q_j,p^*(\lambda_i)\bullet,  0)=0$$ holds in $\cH\otimes[\{q_j\}]$.
\end{theorem}

{\it Proof of Theorem \ref{main}} In fact, we only have to prove
that the equation (\ref{relation}) holds, for any $1\le k \le l$:
by a result of Siebert and Tian (see [13]), these generate the
entire ideal of relations. We apply Theorem \ref{polynomial}  for
the polynomial $D_k$. The only thing which still has to be checked
is the degree condition (\ref{deg}). By the definition of $D_k$,
its degree with respect to $
 Q_j^{1/2}, \Lambda_i$
is the same as the degree of $\Omega_k$ hence, by Theorem
\ref{thmgw2} (ii),  it is equal to the degree of $u_k$. Now, the
inequality
$$\label{ineq}{\rm deg}(u_k)\leq 2\sum_{i=1}^lm_i$$
can be easily verified by checking usual tables: see for instance
\cite{hu2} for max(deg$(u_k)$),  the maximal degree of the
fundamental generators,  and \cite{gw1} for the coefficients $m_i$
from
$$-\alpha_{l+1}^{\vee}= \sum_{i=1}^l m_i\alpha_i^{\vee}.$$ We
summarize the results in the following table. Note that the
inequality (\ref{ineq}) is true for any type of $K$, even though
we need it only for the types $A$, $B$ and $C$ (see Remark 3 in
the introduction).

\vspace{0.5cm}

\begin{tabular}{|c|c|c|}
\hline
Type & max(deg$(u_k)$) & $\sum m_i$\\
\hline
$A_l$ & $l+1$ & $l$\\
$B_l$ & $2l$ & $2l-2$\\
$C_l$ & $2l$ & $l$\\
$D_l$ & $2l-2$ & $2l-3$\\
$E_6$ & $12$ & $11$\\
$E_7$ & $18$ & $17$\\
$E_8$ & $30$ & $29$\\
$F_4$ & $12$ & $8$\\
$G_2$ & $6$  & $5$\\
\hline
\end{tabular}

\hfill$\square$

\section{ The proof of Theorem \ref{polynomial}}

In this section, $q_j$ is always regarded as $e^{t_j}$, $1\le j\le
l+1$. The basic ingredient of the proof are the endomorphisms
$p^*(\lambda_i)\bullet$, $1\leq i\leq l$, of $H^*(LK/T)$ and their
duals with respect to the intersection pairing
$$\langle  \sigma_w,\bar{C}_v\rangle =\delta_{vw},\qquad
v,w\in \Wa,$$
 (see section 4). More precisely, we will denote
by $A_i$ the endomorphism of $H_*(LK/T)$ defined by
$$\langle p^*(\lambda_i)\bullet \sigma_w,
\bar{C}_v \rangle =\langle \sigma_w, A_i(\bar{C}_v)\rangle,
\qquad v,w\in \Wa.$$
 Let us consider an ordering of $\Wa$ compatible with the length;
this will induce an ordering of the basis $\{\sigma_w ~|~ w\in
W_a\}$ of $H^*(LK/T)$ and also of the basis $\{[\bar{C}_w] ~|~
w\in \Wa\}$ of $H_*(LK/T)$. The matrix of $A_i$ with respect to
the latter basis is the transpose of the matrix of
$p^*(\lambda_i)\bullet$ with respect to the former one.

The basis described above induces a natural identification of $H_*(LK/T)$ with the space
$\bR^{\infty}_0$ of all sequences with finite support.
The endomorphisms  of $H_*(LK/T)$ will be represented by matrices of the
form
$$\left (
\begin{matrix}
 a_{11} &  a_{12}  &  \ldots  \\
a_{21}  &  a_{22}  &  \ldots  \\
\vdots    &  \vdots    &\ddots \\
\end{matrix}
\right ),$$
where  each column is in $\bR^{\infty}_0$. Let us denote the space
of  these matrices by $M^{\infty\times \infty}_0(\bR)$.

The matrices $A_i$, $1\leq i \leq l$ have the following features:
\begin{itemize}
\item[a)]  $A_i$ commutes with $A_j$, for any two $i,j$;
\item[b)] $\partial_iA_j=\partial_jA_i$, $1\leq i,j\leq l$;
\item[c)] for any $i\in\{1,\ldots,l\}$, we have that
$$A_i=A'_i(e^{t_j})+A''_i$$
where $A_i'$ is strictly lower triangular, and its coefficients
are linear combinations of $e^{td}=e^{t_1d_1} \ldots e^{t_{l+1}
d_{l+1}}$, $d=(d_1 , \ldots , d_{l+1})\geq 0$, where
$$d_i-m_id_{l+1}\neq 0;$$
 $A''_i$ is strictly upper triangular and its coefficients do not depend
 on $t$
(in particular, the  diagonal of $A_i$ is identically zero).
\end{itemize}

Property a) follows from the associativity of $\bullet$. As for b)
and c), they are direct consequences of the requirements (vi) and
(vii) from Theorem \ref{main}.

The following  result will be needed later:

\begin{proposition}\label{prop} Let $\mathcal{R}$ be a commutative, associative real
algebra with unit and  $A_1, \ldots , A_l\in M^{\infty\times \infty}_0(\cR[e^{t_j}])$ be
 matrices which satisfy the properties
a), b), c) from above. Consider the following system of differential equations:
$$(*) \ \partial_i g=A_ig, \ 1\leq i \leq l$$
where
 $$g=\sum_{d=(d_1, \ldots , d_{l+1})\geq 0}g_d(t_1, \ldots , t_l)e^{td}
\in \cR^{\infty }_0[\{t_i\}][[\{e^{t_j}\}]].$$
The solutions of (*) are uniquely determined by those $g_d^0$ (the degree zero term of
the polynomial $g_d$)
where $d$ satisfies
$$d_1=m_1d_{l+1}, \ldots , d_l=m_ld_{l+1}.$$
More precisely, the solution of (*) can be obtained by the
following recurrence relations: \begin{equation}\label{rec1} g_d =
\left\{
\begin{array}{ll}
                              G_d(g_d^0,g_{d'}|d'<d), \ {\rm if} \  \ed  \\
                              G_d(g_{d'}|d'<d), \ {\rm if \ contrary}
                             \end{array}
                      \right. \end{equation}
where the maps $G_d$ are ``$\cR$-linear in coefficients'' ( by
this I  mean
 that the
coefficients of the polynomial $g_d$ are of the form $C_0g_d^0
+\sum_{j}C_jv_j$ or $\sum_jC_jv_j$ where $C_j\in M^{\infty \times
\infty}_0(\cR)$ and  $v_j\in \cR^{\infty}_0$ are coefficients of
certain polynomials $g_{d'}$ with $d'<d$ (i.e. $d'\leq d$ and
$d'\neq d$).
\end{proposition}

In order to prove Proposition \ref{prop}, we will need the
following lemma:

\begin{lemma}\label{l=1}
Let  $A\in M^{\infty \times \infty}_0(\cR)$ be a matrix of the type
described above and
$g\in\cR ^{\infty}_0[t]$ a polynomial. Consider the differential equation:
$$\frac{{ d} f}{{ d} t}=Af +g,$$
where $f$ is in $\cR ^{\infty}_0[t]$.

(i) If $A-I$ is strictly upper triangular, then we have a unique solution, which depends
$\cR$-linearly on the coefficients of $g$.

(ii) If $A$ is strictly upper triangular, then the solution is uniquely
determined by the degree zero term $f_0\in \cR^{\infty}_0$ of $f$. In fact, the solution
depends $\cR$-linearly on $f_0$ and  the coefficients of $g$.
\end{lemma}

\begin{proof} Put $g=\sum_{k=0}^p g_k t^k$. Look for $f$ of the form $\sum_{k=0}^m f_k t^k$.
We must have $$ f_1 +2 f_2t +\ldots + mf_mt^{m-1}=
  Af_0 +Af_1 t +\ldots + Af_mt^m +g_0 +g_1t +\ldots +g_pt^p,
$$
 hence:
\begin{align}\label{rec2}{}&f_1=Af_0+g_0\nonumber\\
{}& 2f_2=Af_1+g_1=A^2f_0+Ag_0+g_1\nonumber\\
{}& \vdots \\
{}&mf_m=Af_{m-1}+g_{m-1}=c_mA^mf_0+c_{m}A^{m-1}g_0 +\ldots +c_2 A
g_{m-2} +g_{m-1},\nonumber \end{align}
 where $c_j$ are certain nonzero coefficients which
 arise in an obvious way and  $g_r:=0$ for $r\geq p$
(consequently, we have that  $m\geq p$).

If $A$ is the sum of $I$ with a strictly upper triangular matrix, then
 $m$ equals  $p$.
The coefficients $f_j$ are uniquely determined by the recursion
formulae (\ref{rec2}) and the requirement
$$Af_p+g_p=0,$$
which determines the initial term $f_0$.

If $A$ is strictly upper triangular, then we consider the following relations:
$$ (k+p+1)f_{k+p+1}=c_{k+p+1}A^{k+p+1}f_0+ c_{k+p+1}A^{k+p}g_0 +\ldots +
 c_{k+1} A^{k}  g_p,\quad k\geq 0. $$
 If $k\geq 0$ is minimal with the property that
$c_{k+p+1}A^{k+p+1}f_0+ c_{k+p+1}A^{k+p}g_0 +\ldots + c_{k+1}
A^{k}  g_p=0$, then $m$ must be $k+p$ and the coefficients $f_j$
are uniquely determined by $f_0$.

\end{proof}
 {\it Proof of Proposition \ref{prop}} We will   prove this result by
 induction on $l\geq 1$.
Take first $l=1$. We must solve the equation
$$\partial_1g=A_1g,$$
where $g=\sum_{d=(d_1,d_2)\geq 0}g_de^{td}$, $g_d\in \cR^{\infty}_0[t_1]$.
Decompose $A_1$ as $\sum_{d\geq 0}(A_1)^de^{td}$, where $(A_1)^d\in M^{\infty\times \infty}_0(\cR)$.
We identify the coefficients of $e^{td}$ and write the differential equation piecewise as:
$$\frac{{ d} g_d}{{ d} t_1}+(d_1-m_1d_2)g_d=\sum_{d'+d''=d} (A_1)^{d''}g_{d'}.$$
We apply Lemma \ref{l=1}, bearing in mind that $(A_1)^0$ is
strictly upper triangular.

The induction step from $l-1$ to $l$ now follows. We write the $g$
we are looking for as
$$g=\sum_{d=(d_1, \ldots , d_{l+1})\geq 0}g_de^{t_1d_1 +\ldots + t_{l+1}d_{l+1}}=\newline
\sum_{r=(r_1, \ldots ,r_{l-1},r_{l+1})\geq 0}h_re^{t_1r_1+\ldots +t_{l-1}r_{l-1}+
t_{l+1}r_{l+1}},$$
where $g_d\in\cR^{\infty}_0[t_1,\ldots ,t_l]$ and $h_r\in\cR^{\infty}_0[t_l][[e^{t_l}]][t_1, \ldots ,t_{l-1}]$.

Recall that the relation between $g_d$ and $h_r$ is
$$h_r=\sum_kg_{r_1,\ldots ,r_{l-1},k,r_{l+1}}e^{kt_l}.$$
A more detailed formula will give the coefficient
 $h_r^{n_1,\ldots ,n_{l-1}}$
of $t_1^{n_1}\ldots t_{l-1}^{n_{l-1}}$ in $h_r$ as:
\begin{equation}\label{sum}h_r^{n_1,\ldots ,n_{l-1}}=\sum_k \sum_m g_{r_1,\ldots
,r_{l-1},k,r_{l+1}}^{n_1,\ldots ,n_{l-1},m}t_l^me^{kt_l}.\end{equation}
If $h_r^{n_1,\ldots ,n_{l-1}}$ are in $(\cR[t_l][[e^{t_l}]])^{\infty}_0$, then all $g_{r_1,\ldots
,r_{l-1},k,r_{l+1}}$ are in $\cR_0^{\infty}$ (but not conversely).

Now we put $\cS=\cR[t_l][[e^{t_l}]]$ and solve the   system of
differential equations
$$\partial_ih=A_ih, \ 1\leq i \leq l-1,$$
where $h\in \cS_0^{\infty}[t_1,\ldots ,t_{l-1}][[e^{t_1},\ldots
,e^{t_{l-1}},e^{t_{l+1}}]]$. If we regard $A_1,\ldots ,A_{l-1}$ as
elements of $M^{\infty\times \infty}_0(\cS[e^{t_1},\ldots ,
e^{t_{l-1}}, e^{t_l}])$, we can easily see that the conditions a),
b) and c) are still satisfied. By the induction hypothesis, we
know that the latter system can be solved by formulae of the
following type:
$$ h_r =  \left\{ \begin{array}{ll}
                              H_r(h_r^0,h_{r'}|r'<r), \ {\rm if} \  \er  \\
                              H_r(h_{r'}|r'<r), \ {\rm otherwise.}
                             \end{array}
                      \right. $$
We take first the situation when at least one of the conditions
$$\er$$ is not satisfied and show how can we obtain all polynomials
of the form $g_{r_1,\ldots ,r_{l-1},k,r_{l+1}}$, $k\geq 0$ from
$$h_r=H_r(h_{r'}|r'<r).$$
We only have to take into account  (\ref{sum}) and the fact that $H_r$ is a $\cS=
\cR[t_l][[e^{t_l}]]$-linear function  of the coefficients of $h_{r'}$, $r'<r$.
From
$$h_r^{n_1,\ldots ,n_{l-1}}=\sum_{r'<r}B_{r'}h_{r'}^{n'},$$
with $B_{r'} \in M_0^{\infty}(\cR[t_l][[e^{t_l}]])$ we  deduce
easily, by equating the coefficients of $t_l^me^{kt_l}$, that
$g_{r_1,\ldots ,r_{l-1},k,r_{l+1}}^{n_1,\ldots ,n_{l-1},m}$ is a
$\cR$-linear function of the coefficients of $g_{r_1',\ldots,
r_{l-1}',k', r_{l+1}'}$, with $r'<r$ and $k'\leq k$.

If $\er$, we choose $h_r^0\in \cS^{\infty}_0$ to be a solution of:
\begin{equation}\label{sides}\partial_l(h_r^0e^{t_1r_1+\ldots +t_{l-1}r_{l-1}+t_{l+1}r_{l+1}})=
\sum_{r'+r''=r}(A_l^{r'}h_{r''}^0)e^{t_1r_1+\ldots +t_{l-1}r_{l-1}+t_{l+1}r_{l+1}},
\end{equation}
where the matrices $A_l^r\in M^{\infty\times
\infty}_0(\cR[e^{t_l}])$  come from the decomposition
$$A_l=\sum_{r=(r_1,\ldots ,r_{l-1},r_{l+1})\geq 0}A_l^re^{t_1r_1+\ldots +t_{l-1}r_{l-1}+
t_{l+1}r_{l+1}}.$$

In order to be more precise, we write
$$h_r^0=\sum_{k\geq 0}f_k(t_l)e^{kt_l}$$
and then we identify the coefficients of $e^{t_1r_1+\ldots +t_{l-1}r_{l-1}+t_lk+t_{l+1}r_{l+1}}$
in both sides of (\ref{sides}).
One obtains the following sequence of differential equations:
\begin{equation}\label{de}\frac{{ d}f_k}{{d}t_l}+(k-m_lr_{l+1})f_k=(A_l^0h_r^0)_k+b=
\sum_{u+v=k}(A_l^0)^uf_v+b
\end{equation}
where $(A_l^0h_r^0)_k$ denotes the coefficient of $e^{kt_l}$ in
$A_l^0h_r^0$ and $(A_l^0)^u\in M^{\infty\times \infty}_0(\cR)$ is
the coefficient of $e^{t_lu}$ in $A_l^0\in M_0^{\infty\times
\infty}(\cR[e^{t_l}])$. Here $b\in \cR_0^{\infty}[t_l]$ is
obtained by adding the coefficients of $e^{t_lk}$ from all
expressions $A_l^{r'}h_{r''}^0$ where $r'+r''=r$, $r'\neq 0$;
hence it
 depends $\cR$-linearly on the polynomials $g_{r_1'',\ldots ,r_{l-1}'',k',r_{l+1}''}$
where $r''< r$ and $k'\leq k$.

We solve the sequence of differential equations (\ref{de}) using
Lemma \ref{l=1}. First we write (\ref{de}) as:
$$\frac{{ d}f_k}{{ d}t_l}+(k-m_lr_{l+1})f_k=(A_l^0)^0f_k+c+b,$$
where $c\in \cR^{\infty}_0[t_l]$ depends $\cR$-linearly on $f_0,\ldots ,f_{k-1}$.
The matrix $(A_l^0)^0\in M^{\infty\times\infty}_0(\cR)$ is obviously $A_l''$ (see condition c)), hence it is strictly
upper triangular.

We can easily see  that the coefficients of $f_k$  depend $\cR$-linearly on:
\begin{itemize}
\item the coefficients of
$g_{r_1',\ldots ,r_{l-1}',k',r_{l+1}'}$, with $r'<r$ and $k'\leq k$,
if $k\neq m_lr_{l+1}$;

\item  $g_{r_1, \ldots ,r_{l-1},k,r_{l+1}}^0$ and  the coefficients of
$g_{r_1',\ldots ,r_{l-1}',k',r_{l+1}'}$, with $r'<r$ and $k'\leq k$,
if $k=m_lr_{l+1}$.
\end{itemize}

As in the previous case, we take
$$h_r=H_r(h_r^0,h_{r'}|r'<r)$$
and obtain expressions for the coefficients of  $g_{r_1,\ldots
,r_{l-1},k,r_{l+1}}$ via (\ref{sum}) and the structure of $f_k$
described above. More precisely, we take the equality
$$h_r^{n_1,\ldots ,n_{l-1}}=C_0h_r^0+\sum_{r'<r}B_{r'}h_{r'}^{n'},$$
with $C_0, B_{r'} \in M^{\infty\times\infty}_0(\cR[t_l][[e^{t_l}]])$
 and identify the coefficients
of $t_l^me^{t_lk}$.
We deduce  that $g_{r_1,\ldots ,r_{l-1},k,r_{l+1}}^{n_1,\ldots ,n_{l-1},m}$
is
\begin{itemize}
\item  a $\cR$-linear function of the coefficients of $g_{r_1',\ldots, r_{l-1}',k',
r_{l+1}'}$, with $r'<r$ and $k'\leq k$, if $k\neq m_lr_{l+1}$;
\item  a $\cR$-linear function of $g_{r_1, \ldots ,r_{l-1},k,r_{l+1}}^0$ and
 the coefficients of $g_{r_1',\ldots, r_{l-1}',k',
r_{l+1}'}$, with $r'<r$ and $k'\leq k$, if $k= m_lr_{l+1}.$
\end{itemize}

The only thing which remains to be  proved is that the $g$  we
have just
 constructed satisfies
\begin{equation}\label{last}\partial_l g=A_lg.
\end{equation}
To this end, we notice first that
$$\partial_i (\partial_l g-A_lg)=A_i(\partial_l g-A_lg),$$ for all $1\leq i \leq l-1$.
Also $\partial_l g-A_lg$ can be written as
 $$\partial_l g-A_lg=\sum_{r=(r_1,\ldots ,r_{l-1},r_{l+1})\geq 0}q_re^{t_1r_1+\ldots +t_{l-1}r_{l-1}+t_{l+1}r_{l+1}},$$
where $q_r\in \cS^{\infty}_0[t_1, \ldots , t_{l-1}]$, $r\geq 0$
(we have in mind that $\partial_l$ maps \newline
$\cS^0_{\infty}[t_1, \ldots, t_{l-1}][[e^{t_1}, \ldots,
e^{t_{l-1}},e^{t_{l+1}}]]$ onto itself). The degree zero term
$q_r^0$ of $q_r$ is the coefficient of $e^{t_1r_1+\ldots
+t_{l-1}r_{l-1}+t_{l+1}r_{l+1}}$ in
$$\partial_l (h_r^0 e^{t_1r_1+\ldots +t_{l-1}r_{l-1}+t_{l+1}r_{l+1}})-\sum_{r'+r''=r}A_l^{r'}h_{r''}^0e^{t_1r_1+\ldots +t_{l-1}r_{l-1}+t_{l+1}r_{l+1}}.$$
Whenever $\er$ happens, from the choice of $h_r^0$ it follows that
$q_r^0=0$. But these $q_r^0$'s determine $\partial_lg-A_lg$
uniquely, hence the latter  is zero. \hfill $\square$

We  use Proposition \ref{prop} with  $A_i$ replaced by $\frac{1}{h}A_i$ and obtain:

\begin{corollary}
For any $w\in \Wa$ one can find at least one
$s_w\in \cH\otimes \bR[\{t_i\}][[\{e^{t_j}\}]]$ with the properties
$$ \left\{ \begin{array}{ll}
                             h \partial_i s_w =A_i s_w, \ 1\leq i \leq l \\
                             (s_w)_0^0=[\bar{C}_w].
                              \end{array}
                      \right. $$
\end{corollary}

The following lemma will make clear the importance of the formal
series $s_w$, $w\in \Wa$, in proving Theorem \ref{polynomial}.

\begin{lemma}
(i) Let $D\in\bR[ Q_j,\Lambda_i, h]$. Suppose that there exists
$M$ a positive integer with
$$\deg _{Q_j^{1/2},\Lambda_i}(D)\leq 2M$$ and
$$(D(e^{t_j}, h\partial_i,h)\langle 1,s_w \rangle )_d=0,$$
 for any $w\in \Wa$, any multi-index $d$ with $|d|\leq M$ and any
 $h \neq 0.$
Then we have that $$D( e^{t_1}, \ldots ,
e^{t_{l+1}},p^*(\lambda_1)\bullet, \ldots ,p^*(\lambda_l)\bullet,
0)=0.
$$

(ii) We have that $${H}(e^{t_j}, h\partial_i,h) \langle 1,
s_w\rangle=0.$$
\end{lemma}

\begin{proof}
(i) For any $f \in \cH \otimes \bR[e^{t_1}, \ldots , e^{t_{l+1}}]$
 we have that:
$$\begin{array}{lll}h\partial_i\langle f, s_w\rangle&=
\langle h\partial_i f, s_w\rangle +\langle f, h \partial_i
s_w\rangle= \langle h\partial_i f, s_w\rangle +\langle f, A_i
s_w\rangle \\{}&=\langle h\partial_i f, s_w\rangle +\langle
p^*(\lambda_i) \bullet f,  s_w\rangle= \langle (h\partial_i
+p^*(\lambda_i)\bullet) f, s_w\rangle
\end{array}$$
We deduce that
\begin{equation}\label{four}
 D(e^{t_j}, h\partial_i,h) \langle f,s_w \rangle=\langle D(e^{t_1},
 \ldots ,e^{t_{l+1}},p^*(\lambda_1)\bullet+
h\partial_1, \ldots ,p^*(\lambda_l)\bullet+h\partial_l,
h)f,s_w\rangle.\end{equation} Replacing $f$ by $1$ and denoting
${\mathcal D}:=D(e^{t_j},p^*(\lambda_i) \bullet +h
\partial_i,
 h)1$, we obtain:
\begin{equation}
\langle {\mathcal D}, s_w\rangle_d=0,\qquad {\rm if}\ |d|\leq M.
\end{equation}

Notice now that
$${\deg}_{e^{t_j}}{\mathcal D}\leq {\rm deg}_{e^{t_j}}D(e^{t_j},p^*(\lambda_i)\bullet, 0)\leq M.$$
To justify the last inequality,  notice that if
$Q^v\Lambda^u$ is a monomial
from $D(Q_j,\Lambda_i, 0)$, with $u\in \bN^l, v\in \bN^{l+1}$, then we have:
$${\rm deg}_{e^{t_i}}\lambda^{\bullet u}e^{tv}\leq\frac{1}{2}|u|+|v|=
\frac{1}{2}(|u|+2|v|)=\frac{1}{2}{\rm deg}_{ Q_j^{1/2},\Lambda_i
}(Q^v\Lambda^u) \leq \frac{1}{2}\cdot 2M=M.$$

For the rest of the proof, ``degree" and ``degree zero term" will refer  to
the ``variables"  $e^{t_1}, \ldots , e^{t_{l+1}}$.
Decompose ${\mathcal D}$
as
$${\mathcal D}={\mathcal D}_0+{\mathcal D}_1+\ldots +{\mathcal D}_m$$
where ${\mathcal D}_k\in \cH \otimes \bR [e^{t_1},\ldots ,
e^{t_{l+1}}]$
 denotes the sum of all terms of degree $k$,
$0\leq k \leq m$. Recall that the degree zero term of $s_w$ is
$(s_w)_0\in H_*(LK/T)\otimes \bR[t_1, \ldots , t_l]$ with
$(s_w)_0^0=[\bar{C}_w]$. The degree zero term of $\langle
{\mathcal D}, s_w\rangle$ is $\langle {\mathcal D}_0,
(s_w)_0\rangle$. From the vanishing of the latter we obtain that
$$\langle {\mathcal D}_0, (s_w)_0^0\rangle=\langle {\mathcal D}_0,
[\overline{C}_w]\rangle =0,$$ for all $w\in \Wa$, hence ${\mathcal
D}_0=0$.

Also the sum of the terms of degree 1 in
$\langle {\mathcal D}, s_w\rangle$ is 0, hence we have
$$\langle {\mathcal D}_1, (s_w)_0\rangle =0.$$ Exactly as before,
this implies ${\mathcal D}_1=0$. Since $m\leq M$, we can continue
the algorithm until we obtain $\cD_m=0$, hence
$$\cD=0.$$
Now we  let $h$ vary and  deduce
that
$$D(e^{t_j},p^*(\lambda_i)\bullet,0)=0.$$

(ii) When computing $H(e^{t_j}, h\partial_i,h)\langle
1,s_w\rangle$, only the following two relations will be used:
$$h\partial_i\langle 1, s_w\rangle
=\langle p^*(\lambda_i) ,  s_w\rangle $$
and
$$h^2\partial_j \partial_i\langle 1, s_w\rangle=
\langle p^*(\lambda_j)\bullet p^*(\lambda_i), s_w \rangle.$$
  Hence we have to show that:
$$\sum_{i,j=1}^l\langle \alpha_i^{\vee}, \alpha_j^{\vee} \rangle
p^*(\lambda _i)\bullet p^*(\lambda_j) - \sum_{j=1}^{l+1}\langle
\alpha_j^{\vee}, \alpha_j^{\vee} \rangle q_j =0.$$ This follows
immediately from property (v), Theorem \ref{main} and the fact
that
$$\sum_{i,j=1}^l \langle \alpha_i^{\vee}, \alpha_j^{\vee} \rangle m_im_j=
\sum_{i,j=1}^l \langle m_i\alpha_i^{\vee}, m_j\alpha_j^{\vee} \rangle =
\langle \alpha_{l+1}^{\vee}, \alpha_{l+1}^{\vee} \rangle.$$

\end{proof}

Another important step will be made by the following version of
 Kim's Lemma (see \cite{kim} or \cite{gi}):

\begin{lemma}
Let
$$g=g_0 +\sum_{d > 0} g_d e^{td}\in \bR[\{t_i\}][[\{e^{t_j}\}]]$$ be a formal series
which satisfies
$$g_0=0 \ \ {\rm and} \ \ H(e^{t_j}, h\partial_i,h) g=0.$$
Then $g_d=0$ for any $ d$ with $|d|\leq \sum_{i=1}^lm_i$.
\end{lemma}

\begin{proof}
Take $d\in \bN^{l+1}$ with $g_d\neq 0$ and $|d|:=\sum_{i=1}^{l+1} d_i>0$ minimal.
 From $H(e^{t_j}, h\partial_i,h) g=0$ it follows that
$$(\sum_{i,j=1}^l\langle \alpha_i^{\vee}, \alpha_j^{\vee}\rangle
\partial ^2_{ij})(g_de^{td})=0$$
and then
$$(\sum_{i,j=1}^l\langle \alpha_i^{\vee}, \alpha_j^{\vee}\rangle
\partial ^2_{ij})(e^{td})=0.$$

 On the other hand, we have that
$ (\sum_{i,j=1}^l\langle \alpha_i^{\vee}, \alpha_j^{\vee}\rangle
\partial ^2_{ij})(e^{td})     =
\sum_{i,j=1}^l\langle \alpha_i^{\vee},\alpha_j^{\vee}\rangle (d_i-m_id_{l+1})( d_j-m_jd_{l+1}) e^{td}=
 || \sum_{j=1}^l (d_j-m_jd_{l+1}) \alpha _j^{\vee})||^2 e^{td},$
hence we must have $d_i=m_id_{l+1},$ for all $1\leq i\leq l$. From
the fact that $d\neq 0$, it follows that $d_{l+1}\geq 1$ and
$d_i\geq m_i$ for any $1\leq i \leq l$ and finally that
$$|d|\geq 1+\sum_{i=1}^lm_i.$$
The assertion stated in the lemma follows from the minimality of $|d|$.
\end{proof}

{\it Proof of Theorem \ref{polynomial}} We show that
$$g:=D(e^{t_j}, h\partial_i,h)\langle 1, s_w\rangle =\sum_{d\geq 0}g_de^{td}$$
vanishes. Put first $f=1$ in (\ref{four}) and obtain
\begin{equation}\label{g}\begin{array}{lll}
g&=\langle D(e^{t_1}, \ldots e^{t_{l+1}},p^*(\lambda_1)\bullet+
h\partial_1, \ldots ,p^*(\lambda_l)\bullet+h\partial_l,
h)1,s_w\rangle \\{}&= \langle D (0, \ldots,
0,p^*(\lambda_1)\bullet,\ldots ,p^*(\lambda_l) \bullet, h)+R,
s_w\rangle,\end{array}\end{equation} where $R\equiv 0$ mod
$\{e^{t_j}\}$. Hence $g_0\in \bR[t_1, \ldots , t_l]$ will be the
same as  $$\langle D(0, \ldots, 0,p^*(\lambda_1), \ldots
,p^*(\lambda_l),  h), s_w\rangle_0$$ where, as usual, the
subscript 0 indicates the degree zero term with respect to
$e^{t_1}, \ldots ,e^{t_{l+1}}$. But $$D( 0, \ldots,
0,p^*(\lambda_1), \ldots ,p^*(\lambda_l), h)= p^*(D( 0, \ldots,
0,\lambda_1, \ldots ,\lambda_l, h))=0,$$ hence $g_0=0$.
 If $$M:=\sum_{i=1}^lm_i,$$ then $g_d=0$, for $|d|\leq M$ (we take into
 account that $[D(e^{t_j}, h\partial_i,h), H(e^{t_j}, h\partial_i,h)]=0$,
and also that
 $H(e^{t_j}, h\partial_i,h)\langle 1,s_w\rangle =0$, which imply
 that
$H(e^{t_j}, h\partial_i,h) g =0$ and we apply Lemma 6.5). Finally,
from Lemma 6.4
  it follows that
$$D( e^{t_j},\lambda_i\bullet,0)=0$$
and the proof is complete. \hfill$\square$

\vspace{1cm}

Department of Mathematics, University of Toronto

100 St. George Street, Toronto, Ontario, M5S 3G3 Canada

E-mail address: {\tt amare@math.toronto.edu}

\end{document}